\documentclass[12pt]{article}
\usepackage{latexsym}

\usepackage{graphics,graphicx}
\usepackage{amssymb}
\usepackage{amsmath}
\DeclareGraphicsExtensions{.png}
\usepackage{latexsym}
\usepackage{tikz}
\usepackage{xcolor}
\usepackage{graphics,graphicx}
\usepackage{amssymb}
\usepackage{amsmath}
\DeclareGraphicsExtensions{.png}
\usepackage{latexsym}
\usepackage{tikz}
\usepackage{xcolor}

\newtheorem{theorem}{Theorem}

\newtheorem{lemma}[theorem]{Lemma}
\begin{document}
\title{Optimal Control of Hydro-Electric Power Plants with Uncontrolled Spillways}
\author{Maria do Ros\'{a}rio de Pinho, Maria Margarida A. Ferreira, and Georgi Smirnov 
\thanks{The authors thank the support of  Portuguese Foundation for Science and Technology
(FCT) in the framework of the Strategic Funding UIDB/04650\-/2020,  
 Portuguese Foundation for Science and Technology (FCT/MCTES), its support within the framework of the Associated Laboratory ARISE (LA/P/0112/2020) and the Research and Development Unit SYSTEC through Base (UIDB/00147/2020) and Programmatic (UIDP/00147/2020) funds and the support of  FCT/MCTES (PIDDAC) through project 2022.02801.PTDC-UPWIND-ATOL.\\
 Maria do Ros\'{a}rio de Pinho is with Universidade do Porto, Faculdade de Engenharia, 
DEEC, SYSTEC, Porto, Portugal (e-mail: mrpinho@fe.up.pt). \\
Maria Margarida A. Ferreira was with Universidade do Porto, Faculdade de Engenharia, 
DEEC, SYSTEC, Porto, Portugal (e-mail: mmf@fe.up.pt).\\
Georgi Smirnov is with 
Universidade do Minho, Dep. Matem\'{a}tica, 
Physics Center of Minho and Porto Universities (CF-UM-UP), 
Campus de Gualtar, Braga, Portugal (e-mail: smirnov@math.uminho.pt).}}
\date{~}

\maketitle

\begin{abstract}
In this paper, we study an optimal control problem for a cascade of hydroelectric power plants with reversible turbines and uncontrolled spillways. The system dynamics are governed by a linear control model subject to path constraints. The aim is to maximize the power production profit while respecting operational restrictions on reservoir water levels. The challenge is the presence of uncontrollable spillways: their discontinuous nature and the fact that they are activated at the state boundary prevent the application of known necessary conditions of optimality. To overcome this, we derive necessary conditions by approximating the original problem through a sequence of standard optimal control problems using exponential penalty functions. The applicability of resulting conditions are illustrated by an example.
\end{abstract}

\section{Introduction}
\label{intro}

Cascades of hydro-electric power plants with reversible turbines play an important role in renewable energies and environment protection \cite{Labadie}. Power plants with reversible turbines in multireservoirs systems allow for the production of electric energy during the daytime, while, during the nighttime when the price of the electric energy is low, such systems have the possibility of pumping water from downstream reservoirs using, for example, wind energy.
 
In this paper, we consider a cascade of hydro-electric power plants with reversible turbines and uncontrolled spillways.  In the case of large power plants,  spillways usually are controlled in the sense that, when needed, water is discharged to downstream reservoirs. Controlled spillways are of importance for the optimal managing of power production. In this respect, see, e.g., \cite{MC}, where the reader can find an example of a solution to an optimal control problem with an active spillway. However, in the case of cascades with small power plants, the spillways may be  uncontrollable, like ogee or bell-mouth spillways. This means that spillways may be active only when the water level in the reservoir attains its maximum value.

The profit of power production for hydro-electric plants has been amply studied  in the literature often using optimization techniques.
In the context of optimal control for cascades of hydro-electric power plants, the literature seems to be surprisingly sparse; for example, sufficient conditions are studied in \cite{JOTA1} and some numerical methods can be found in \cite{MC} and  \cite{Bush}. 
 
Here, our aim is to optimize the  profit of power production for cascades of hydro-electric power plants with reversible turbines and uncontrolled spillways.   We  consider this problem in the framework of optimal control theory.  The model we use is based on \cite{VALORAGUA,ruey,MC}.  

Optimal control problems for cascades with controlled spillways are  state constrained problems  for which necessary  conditions of optimality are well-known (see, e.g., \cite{V,BV}). However, the situation changes drastically in the presence of uncontrollable spillways. The problem cannot be reformulated to fit within the existing theoretical frameworks. This is because of the discontinuous nature of the spillway and its activation solely on the boundary of state constraints. As a result, the spillway cannot be treated as a state variable. Moreover, its inherent properties prevent it from being interpreted as either a control variable or a time-dependent parameter. It may also seem that our problem may be reduced to a problem involving sweeping processes in vein of those treated in \cite{JOTA2}, but this does not hold. Thus, another approach is called for. To derive necessary conditions of optimality for such problems, we construct a sequence of standard approximating optimal control problems using exponential penalty functions. By applying the well-known Pontryagin Maximum Principle to this sequence and passing to the limit, we obtain the desired necessary conditions.

The paper is organized in the following manner. In Sec. 2, we introduce notation and the statement of the problem. The construction of the sequence of approximating problems is given in Sec. 3. The main result of this work is formulated in Sec. 4, where we  also present an illustrative  example.  Sec. 5 contains a brief conclusion. The proofs are placed in the Appendix.

\section{Notation and problem statement}\label{sec:1}

Before presenting our main result, we start by introducing the notation used throughout this paper and describing the problem under consideration.

\subsection{Notation}

Throughout this  paper we denote  the set of real numbers by $R$ and 
 the usual  $n$-dimensional  space    of  vectors $x=(x_1,\ldots ,x_n)$,  where $x_i\in R$, $ i=1,\ldots,n$, by $R^n$. 
 The inner product of two vectors $x$ and $y$ in $R^n$ is defined by $\langle x,y\rangle =\sum_{i=1}^n x_iy_i$, and the norm  of  a  vector  $x\in  R^n$  is    $|x|=\langle x,x\rangle^{1/2}$.  If $x\in R^n$, then ${\rm diag}(x)$ is used to denote an $n\times n$ diagonal matrix with the diagonal $x_1,\ldots,x_n$.
 
  The space  $L_{\infty}([a,b]; R^p)$ (or simply $L_{\infty}$ when the domains are clearly understood) is the Lebesgue space of essentially bounded functions  $h:[a,b]\to R^p$ and the space $L_2([a,b]; R^p)$ is the space of Lebesgue square integrable functions.   The space of continuous functions is denoted by $C([a,b]; R^p)$ and the space of absolutely continuous functions  by  $AC([a,b]; R^p)$. We say that $h\in BV([a,b], R^p)$ if $h$ is a function of bounded variation continuous from the right, may be except at the point $a$. So, the dual space $(C([a,b], R^p))^*$ can be identified with $BV([a,b], R^p)$.

\subsection{Problem statement}

\begin{figure}[ht!]

\centering

\begin{tikzpicture}

\draw [thick] (-2,2) -- (0,0);

\draw[thick] [<->] (-1.5,1.5) -- (-0.5,0.5);

\draw [thick] (2,2) -- (0,0);

\draw[thick] [<->] (1.5,1.5) -- (0.5,0.5);

\draw [thick] (0,2) -- (0,0);

\draw[thick] [<->] (0,2.9) -- (0,1.1);

\draw [thick] (0,5) -- (0,2);

\draw[thick] [->] (0,5) -- (0,4.5);

\draw [thick] (0,0) -- (0,-2);

\draw[thick] [->] (0,-1) -- (0,-1.9);

\filldraw [black]    (-1,1) circle (5pt);

\node at (-1.2,1) [left] {$2$};

\filldraw [black]    (1,1) circle (5pt);

\node at (1.2,1) [right] {$4$};

\filldraw [black]    (0,2) circle (5pt);

\node at (0.2,2) [right] {$3$};

\filldraw [black]    (0,4) circle (5pt);

\node at (0.2,4) [right] {$1$};

\filldraw [black]    (0,-1) circle (5pt);

\node at (0.2,-1) [right] {$5$};

\end{tikzpicture}

\caption{ Cascade of five hydro-power plants. Three of them, $2$, $3$ and $4$, marked  with double sided arrows, have reversible turbines.}

\label{fig:1}

\end{figure}
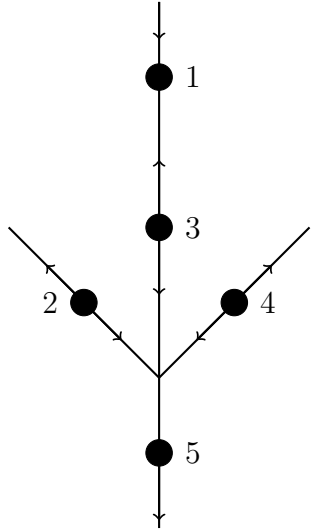

Here, we consider a problem  associated to a system of $\cal I$  hydro-electric power plants with a cascade structure. Figure \ref{fig:1} represents a scheme for such  cascade with  five power plants.

For a  cascade of hydro-electric power plants, we assume that  the dynamics of water volumes, $V_i(t)$, in the reservoir $i=1,2, \ldots, {\cal I}$,  is des\-cri\-bed by the following control system
\begin{equation}
\dot{V}_i(t) = A_i-u_i(t)-s_i(t)+\sum_{j\in J(i)}\left( u_j(t)+s_j(t)\right) ,\label{*}
\end{equation}
where  $J(i)$ represents the set (that may be empty) of reservoir indexes $j<i$ from which the water flows to  reservoir $i$, $u(t)= (u_1(t), \dots, u_{\cal I}(t))$ are the controls  representing the turbined/pumped volumes of water  for each reservoir at time $t$, and $A_{i}$ are the incoming flows, $\ i=1,2,\dots,{\cal I}$. The  spillways   for each reservoir at time $t$ are denoted by $s(t)= (s_1(t), \dots, s_{\cal I}(t))$. The equation \eqref{*} is called \emph{the water balance equation}; it is widely used in the literature (see, e.g.,  \cite{mariano}).

We   assume throughout that $J(i)\cap J(i')=\emptyset$, whenever $i\neq i'$. For example, in the case presented in Figure \ref{fig:1}, we have $J(5)=\{ 2,3,4\}$ and $J(3)=\{ 1\}$. Hence, the index $i=J^{-1}(j)>j$ is associated to a unique downstream reservoir (if it exists) which receives water from reservoir $j$;  in other words, $j\in J(i)$. 

The control variables and the water volumes satisfy the following technical constraints:
\begin{equation}
\label{**}
 V_i(0)=V_i(T),~~ V_i(t)\in [V_i^m,V_i^M], ~~ u_i(t)\in [u_i^m,u_i^M].
\end{equation}
Here $V_i^m$ and $V_i^M$,  $i=\overline{1,{\cal I}}$, stand for the imposed minimum and maximum water volumes;   $u_i^m$ and $u_i^M$,  $i=\overline{1,{\cal I}}$, are the allowed minimum and maximum  turbined/pumped water flows. 

The spillways $s_i(t)$, $i=\overline{1,{\cal I}}$, satisfy the conditions
\begin{equation}
\label{***}
s_i(t)\geq 0,\;\;\; s_i(t)(V_i(t)-V_i^M)=0,
\end{equation}
i.e., $s_i(t)$ can be positive only if $V_i(t)=V_i^M$.

For simplicity, we assume that the reservoirs have cylindrical form and that the gravity constant is equal to one. Also, we assume that  all the potential energy is converted into electric energy.

The objective is to find the optimal controls $\hat{u}_i(\cdot ) \in L_{\infty}([0,T], R)$  and respective volumes $\hat{V}_i(\cdot ) \in AC([0,T], R)$  giving an optimal strategy for the management of water in order to produce the maximum financial return of energy production in the system:
$$
\text{maximize} \; J(u(\cdot ),V(\cdot ))
$$
\begin{equation}
\label{cost}
=\sum_{j=1}^{{\cal I}}\int^{T}_{0}c(t)u_j(t)\left(\frac{V_j(t)}{S_j}+h_j
 -\frac{V_{J^{-1}(j)}(t)}{S_{J^{-1}(j)}}-h_{J^{-1}(j)}\right)dt.
\end{equation}

Here, $c(\cdot )$ is the price of the energy, $h_j$, $j=\overline{1,{\cal I}}$,  denote  the elevations above mean sea level  of the reservoir bases,  and $S_j$, $j=\overline{1,{\cal I}}$, are the areas of the  bases of the reservoirs.  In \eqref{cost}, the price $c(t)$ is multiplied by an expression representing potential energy at time $t$ completely transformed into  electrical energy. Observe that, if $J^{-1}(j)=\emptyset$, then respective terms in \eqref{cost} do not appear.

\vspace{5mm}

Although problem \eqref{*} -- \eqref{cost} may seem to be a standard optimal control problem with state constraints, the presence of uncontrolled spillways, functions to be determined, that may be different from zero only when some state constraints are violated, adds serious difficulty as far as necessary conditions are concerned. This is because  known results do not apply to this problem. To overcome such difficulties we resource to  penalty function methods, a widely applied technique to derive necessary conditions of optimality for optimal control problems (see, e.g., \cite{Pappas,Bayen}, to name but a few). Since spillway variables may only be active on the state constraint boundary, we opt to use exponential penalty function techniques in vein of those previously used for different systems in \cite{SVVA1,JOTA2,SCL}. An easy to follow introduction to this method can be found in the survey \cite{S}, where the reader can find probably the most important references to the problems covered there. 

\vspace{5mm}

Problem \eqref{*} -- \eqref{cost} can be written in vector notation, in the following form:
\begin{eqnarray}
&& \int_0^T c(t)\langle u(t), S^{-1}V(t)+h \nonumber \\
&&-{\cal M}^*(S^{-1}V(t)+h)\rangle dt \rightarrow\max,\label{P1}\\
&&  \dot{V}=A-u-s+{\cal M}(u+s),\label{P2}\\
&&  u\in [u^m,u^M],\label{P3}\\
&& V\in [V^m,V^M],\;\; s\geq 0,\;\; {\rm diag} (s)(V-V^M)=0,\label{P4}\\
&& V(0)=V(T),\label{P5}
\end{eqnarray}
where
\begin{itemize}
\item incoming flows $A=(A_1,\ldots,A_{\cal I})$;
\item minimum volumes $V^m=(V^m_1,\ldots,V^m_{\cal I})$;
\item maximum volumes $ V^M=(V^M_1,\ldots,V^M_{\cal I})$;
\item minimum  turbined/pumped  volumes $$u^m=(u^m_1,\ldots,u^m_{\cal I}); $$
\item maximum  turbined/pumped  volumes $$u^M=(u^M_1,\ldots,u^M_{\cal I});$$
\item elevation of reservoir bases $h=(h_1,\ldots,h_{\cal I})$;
\item areas of of reservoir bases $$S={\rm diag}(S_1,\ldots,S_{\cal I})\in R^{{\cal I}\times {\cal I}},$$
\end{itemize}
with $A$, $V^m$, $V^M$, $u^m$, $u^M$, $h$ and $S$ are in $R^{\cal I}$, 
\begin{itemize}
\item reservoir volumes $$V(\cdot)=(V_1(\cdot),\ldots,V_{\cal I}(\cdot))\in AC([0,T],R^{\cal I});$$
\item turbined/pumped volumes $$u(\cdot)=(u_1(\cdot),\ldots,u_{\cal I}(\cdot))\in L_{\infty}([0,T],R^{\cal I});$$
\item reservoir spillways $$s(\cdot)=(s_1(\cdot),\ldots,s_{\cal I}(\cdot))\in L_{\infty}([0,T],R^{\cal I})$$
\end{itemize}
and matrix ${\cal M}\in R^{{\cal I}\times {\cal I}}$.
Observe that  the unknown spillways $s(\cdot)$ are  essentially bounded measurable functions.

\vspace{5mm} 

The matrix ${\cal M}\in R^{{\cal I}\times {\cal I}}$ has a special structure describing the topology of a cascade of hydro-electric power plants under consideration. Thus,  ${\cal M}$  is such that  ${\cal M}e_{j}=e_{J^{-1}(j)}$, where $e_j$, $j=1,\ldots,{\cal I}$, are the vectors of the orthonormal basis (recall that the reservoir $J^{-1}(j)$ receives water from the reservoir $j$.). 

Take, for example,  the cascade in Figure \ref{fig:1}.
The respective matrix $\cal M$ has the form
$$
{\cal M}=\left(
\begin{array}{ccccc}
0&0&0&0&0\\
0&0&0&0&0\\
1&0&0&0&0\\
0&0&0&0&0\\
0&1&1&1&0
\end{array}
\right).
$$

\section{Preliminary results}

Our aim is to establish necessary conditions of optimality for problem \eqref{P1} -- \eqref{P5}. As mentioned in the Introduction, a sequence of approximating standard optimal control problems plays a central role in our analysis. Before presenting the necessary conditions in the next section, we first construct this sequence. The proposed family of approximating problems, along with their corresponding results, is not only instrumental in deriving the necessary conditions but also of independent interest, for instance, as a basis for developing numerical methods applicable to related problems.

\subsection*{Approximating system}

Let us focus on the equations
\begin{equation}
\label{s1}
 \dot{V}=A-u-s+{\cal M}(u+s)
\end{equation}
 and
\begin{equation}
\label{s2}
V\leq V^M,\;\; s\geq 0,\;\; {\rm diag} (s)(V-V^M)=0.
\end{equation}
  
\vspace{5mm}

 Our starting point is the construction of a sequence of approximating systems.  This is done in the following lemma, where we state the convergence of the solutions to the approximating systems to a solution of problem  \eqref{s1} -- \eqref{s2}.
\begin{lemma}
\label{APP}
Consider a sequence $u^{\gamma}$ converging in the weak-* topology of $L_{\infty}([0,T],R^{\cal I})$ to some control $\tilde{u}\in  [u^m,u^M]$, a sequence $V^{\gamma,0}\leq V^M$ converging to $V^0\in R^{\cal I}$ and the Cauchy problems
\begin{eqnarray}
&&  \dot{V}_i^{\gamma}=A_i-u_i^{\gamma}-s^{\gamma}_i +\sum_{j\in J(i)}\left( u_j^{\gamma}+s^{\gamma}_j\right),\nonumber\\
&& {\rm for}\,\; i=1,\ldots,{\cal I},\label{A1} \\
&& V^{\gamma}(0)= {V}^{\gamma,0},\label{A2}
\end{eqnarray}
where $s^{\gamma}_i=\gamma^i e^{\gamma^i(V_i^{\gamma}-V_i^M)}$, $i=1,\ldots,{\cal I}$. (Here and in what follows $\gamma^i$ stands for $\gamma$ to the $i$-th power.)

When $\gamma\to\infty$, the sequence of control processes $(u^{\gamma},V^{\gamma},s^{\gamma})$ solving  problem \eqref{A1} -- \eqref{A2}  converges to a solution $(\tilde{u},\tilde{V},\tilde{s})\in L_{\infty}([0,T],R^{\cal I})\times C([0,T],R^{\cal I})\times L_{\infty}([0,T],R^{\cal I})$ of system  \eqref{s1} -- \eqref{s2} with the initial condition  $\tilde{V}(0)=V^0$, where $s^{\gamma}$ converges in weak topology of   $L_{2}([0,T],R^{\cal I})$ and $V^{\gamma}$ converges in the strong topology of  $C([0,T],R^{\cal I})$. 

Moreover, the Cauchy problem for system  \eqref{s1} -- \eqref{s2} has a  solution and this solution is unique.
\end{lemma}

\vspace{5mm}

The proof of the lemma is in the Appendix.

\vspace{3mm}

{\bf Remark.} From the differential equation \eqref{A1} for $V^{\gamma}$,  we see that the constructed solution to the Cauchy problem for \eqref{s1} and \eqref{s2} is bounded from below, i.e.,  there exists some $\check{V}\in R^{\cal I}$ such that 
\begin{equation}\label{Vbound}
\check{V}\leq \tilde{V} \leq V^M.
\end{equation}

\subsection*{Approximating problem} 
 
 Based on the previous lemma, we now turn attention to  a family of approximating optimal control problems. 
 
Let $(\hat{u},\hat{V},\hat{s})$ be a global optimal solution to problem \eqref{P1} -- \eqref{P5}, and let $\alpha>0$, $\epsilon>0$ and $\gamma>0$. 
Let us consider the minimization of the functional
\begin{eqnarray}
&& {\cal J}(u,V)= \frac{1}{2}|V(0)-\hat{V}(0)|^2\nonumber\\
&& - \int_0^T\sum_{j=1}^{{\cal I}} c u_j\left(\frac{V_j}{S_j}+h_j-\frac{V_{J^{-1}(j)}}{S_{J^{-1}(j)}}-h_{J^{-1}(j)}\right) dt\nonumber\\
&&
+ \int_0^T\!\left(\sum_{j=1}^{{\cal I}}\frac{e^{-\gamma (V_j-V_j^m+\epsilon)}-1}{\sqrt{\gamma}}+\alpha |u-\hat{u}|^2\right) dt \label{ap1}
\end{eqnarray}
subject to 
\begin{eqnarray}
&&  \dot{V}_i=A_i-u_i-s^{\gamma}_i \nonumber \\
&& +\sum_{j\in J(i)}\left( u_j+s^{\gamma}_j\right),\;\, i=1,\ldots,{\cal I},\label{ap2}\\
&&  u\in [u^m,u^M],\label{ap3}\\
&& \frac{1}{2}|V(0)-V(T)|^2\leq \epsilon, \label{ap4}
\end{eqnarray} 
where, once again, $s^{\gamma}_i=\gamma^i e^{\gamma^i (V_i^{\gamma}-V_i^M)}$, $i=1,\ldots,{\cal I}$.

A comparison between this problem and problem  \eqref{P1} -- \eqref{P5} is called for. Observe that the periodic condition  \eqref{P5} and the state constraint  \eqref{P4} are slightly relaxed. This is important to guarantee the existence of admissible solutions to the above problem. 
While equation \eqref{P5} becomes now \eqref{ap4}, the constraints  \eqref{P4} are handled now by  
$$
\frac{e^{-\gamma (V_j-V_j^m+\epsilon)}-1}{\sqrt{\gamma}}
$$
in \eqref{ap1} for the state constraint $V\geq V^m$, and  $s^{\gamma}$ in \eqref{ap2} for the constraint $V\leq V^M$.

\vspace{5mm}

\begin{lemma}
\label{lem1}
For large $\gamma$,  problem \eqref{ap1} -- \eqref{ap4} has admissible solutions and, as a consequence, there exists a global optimal control process  $(\tilde{u}^{\gamma,\epsilon,\alpha},\tilde{V}^{\gamma,\epsilon,\alpha},\tilde{s}^{\gamma,\epsilon,\alpha})$. 

When $\gamma\to\infty$, the sequence $\tilde{u}^{\gamma,\epsilon,\alpha}$ converges in the weak-* topology of $L_{\infty}([0,T],R^{\cal I})$ to some control $\tilde{u}^{\epsilon,\alpha}\in  [u^m,u^M]$,  the sequence $\tilde{V}^{\gamma,\epsilon,\alpha}$ uniformly converges to $\tilde{V}^{\epsilon,\alpha}$ and  the sequence $\tilde{s}^{\gamma,\epsilon,\alpha}$ converges in the weak topology of $L_2([0,T],R^{\cal I})$ to some  $\tilde{s}^{\epsilon,\alpha}\in L_{\infty}([0,T],R^{\cal I})$.
The process $(\tilde{u}^{\epsilon,\alpha},\tilde{V}^{\epsilon,\alpha},\tilde{s}^{\epsilon,\alpha})$ is an admissible solution to problem $({\cal P}_{\epsilon,\alpha})$ defined as
$$
\!
\begin{array}{lll}
&& {\cal J}(u,V)= \frac{1}{2}|V(0)-\hat{V}(0)|^2\\
&& - \int_0^T\sum_{j=1}^{{\cal I}} c  u_j \left( \frac{V_j}{S_j}+h_j-\frac{V_{J^{-1}(j)}}{S_{J^{-1}(j)}}-h_{J^{-1}(j)}\right) dt\\
&& +\int_0^T\alpha |u-\hat{u}|^2 dt \rightarrow\min,\label{pe1}\\
&&  \dot{V}_i=A_i-u_i-s_i+\sum_{j\in J(i)}( u_j+s_j),\;\; i=1,\ldots,{\cal I},\label{pe2}\\
&&  u_i\in [u_i^m,u_i^M],\label{pe3}\\
&& V_i\in [V_i^m-\epsilon,V_i^M],\label{pe4}\\
&& s_i\geq 0,\quad s_i(V_i-V_i^M)=0,\label{pe5}\\
&& \frac{1}{2}|V(0)-V(T)|^2\leq \epsilon. \label{pe6}
\end{array}
$$ 
\end{lemma}

\vspace{5mm}

The proof of the lemma is in the Appendix.

\vspace{5mm}

Now, we turn to the relaxations introduced in the periodic end-point constraints and the state constraint $V\geq V^m$.  

\vspace{5mm}

\begin{lemma}
\label{lem2}
Let $\epsilon\to 0$. Then, without loss of generality, the sequence of processes $(\tilde{u}^{\epsilon,\alpha},\tilde{V}^{\epsilon,\alpha},\tilde{s}^{\epsilon,\alpha})$ tends to the process $(\hat{u}^{\alpha},\hat{V}^{\alpha},\hat{s}^{\alpha})$, a global optimal solution to problem $({\cal P}_{\alpha})$ defined as
$$
\!
\begin{array}{lll}
&& {\cal J}(u,V)= \frac{1}{2}|V(0)-\hat{V}(0)|^2\\
&& - \int_0^T\sum_{j=1}^{{\cal I}} c  u_j \left( \frac{V_j}{S_j}+h_j-\frac{V_{J^{-1}(j)}}{S_{J^{-1}(j)}}-h_{J^{-1}(j)}\right) dt\\
&& +\int_0^T\alpha |u-\hat{u}|^2 dt \rightarrow\min,\label{pe1a}\\
&&  \dot{V}_i=A_i-u_i-s_i+\sum_{j\in J(i)}( u_j+s_j),\;\; i=1,\ldots,{\cal I},\label{pe2a}\\
&&  u_i\in [u_i^m,u_i^M],\label{pe3a}\\
&& V_i\in [V_i^m,V_i^M],\label{pe4a}\\
&& s_i\geq 0,\quad s_i(V_i-V_i^M)=0,\label{pe5a}\\
&& V(0)=V(T), \label{pe6a}
\end{array}
$$ 
 where the convergence is understood in the following sense: the sequence $(\tilde{u}^{\epsilon,\alpha},\tilde{s}^{\epsilon,\alpha})$ converges in weak-* topology of $L_{\infty}([0,T],R^{\cal I}\times R^{\cal I})$ and the sequence $\tilde{V}^{\epsilon,\alpha}$ converges in the strong topology of $C([0,T],R^{\cal I}\times R^{\cal I})$.
\end{lemma}

\vspace{5mm}

The proof of the lemma is in the Appendix.

\section{Necessary conditions of optimality}

The main result of this work is contained in the following theorem establishing necessary  conditions of optimality for the problem under consideration. 

\vspace{5mm}

\begin{theorem}
\label{th_main}
Let $(\hat{u},\hat{V},\hat{s})$ be a global optimal solution to problem \eqref{P1} -- \eqref{P5}. Then there exist $\lambda\geq 0$, $(\mu,\xi,p)\in BV([0,T],R^{\cal I}\times R^{\cal I}\times R^{\cal I})$ satisfying the following relations:
$$
d{p}_i=d\xi_i-\lambda c \frac{\hat{u}_i-\sum_{j\in J(i)}\hat{u}_j}{S_i} dt- d\mu_i,
$$
$$
 (\hat{V}_i-V_i^M)d\xi_i=0,\; (\hat{V}_i-V_i^m)d\mu_i=0,\; d\mu_i\geq 0,\; \mu_i(0)=0,
$$
$$
p(0)=p(T),
$$
$$
\hat{s}_i(p_i-p_{J^{-1}(i)})=0,
$$
$$
\lambda+|p(T)|+\sum_{i=1}^{\cal I} \mu_i(T)=1.
$$
Moreover, almost everywhere the optimal control $\hat{u}(t)$ solves the maximization problem
$$
\sum_{i=1}^{\cal I} p_i\left( -u_i+\sum_{j\in J(i)}u_j \right)
$$
$$
+\lambda  \sum_{j=1}^{\cal I} cu_j\left( \frac{\hat{V}_j}{S_j}+h_j - \frac{\hat{V}_{J^{-1}(j)}}{S_{J^{-1}(j)}}-h_{J^{-1}(j)} \right)\to\max_{u\in U}.
$$ 
\end{theorem}

\vspace{5mm}

The proof of this theorem, which relies on the construction from the previous section, is provided in the Appendix.

An example is presented next to demonstrate the utility of the preceding theorem in identifying the solution of a problem. 

\subsection*{Example}

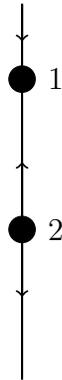
\begin{figure}[ht!]
\centering

\begin{tikzpicture}

\draw [thick] (0,2) -- (0,0);
\draw[thick] [<->] (0,2.9) -- (0,1.1);
\draw [thick] (0,5) -- (0,2);
\draw[thick] [->] (0,5) -- (0,4.5);
\filldraw [black]    (0,2) circle (5pt);
\node at (0.2,2) [right] {$2$};
\filldraw [black]    (0,4) circle (5pt);
\node at (0.2,4) [right] {$1$};

\end{tikzpicture}

\caption{ Cascade of two hydro-power stations. One of them, $2$,  marked with double sided arrows, has reversible turbine.}
\label{fig:2}
\end{figure}

Consider a  cascade of two hydro-power plants, one of them, $2$, with reversible turbine (see Fig \ref{fig:2}). To simplify we assume that $S_1=S_2=1$.    
The problem we consider  is then the following:
\begin{eqnarray*}
&&\int_0^{16} c(u_1(V_1+13-V_2)+u_2V_2)dt\rightarrow\max,\\
&&\dot{V}_1=2-u_1-s_1,\\
&&\dot{V}_2=-u_2-s_2+u_1+s_1,\\
&& V_1(0)=V_1(16),\; V_2(0)=V_2(16),\\
&&u_1\in [0,1],\; u_2\in [-1,3],\; V_1\in [1,5],\; V_2\in [3,12],\\
&& s_1\geq 0,\; s_2\geq 0,\; s_1(V_1-5)=0,\; s_2(V_2-12)=0,
\end{eqnarray*}
where 
$$
c(t)=\left\{
\begin{array}{rl}
3,& t\in [0,6[,\\
11, & t\in [6,16[,\\
3,& t=16.
\end{array}
\right.
$$

In this case, we have $J(2)=\{ 1\}$, $J(1)=\emptyset$ and $J^{-1}(1)=\{ 2\}$. 

The parameters used in the example are purely made up and chosen to guarantee that the first power station has a permanent spillway. This is a situation that we may have during a very rainy  period.
 
 Let $(u_i,V_i,s_i)$, $i=1,2$, be a global optimal solution to this problem. 
Theorem \ref{th_main} asserts the existence of  $\lambda\geq 0$,   $(\mu_i,\xi_i,p_i)\in BV([0,T],R)$, $i=1,2$, such that 
\begin{eqnarray*}
&& dp_1=d\xi_1-\lambda cu_1dt-d\mu_1,\\
&& dp_2=d\xi_2-\lambda c(u_2-u_1)dt-d\mu_2,\\
&& (V_1-5)d\xi_1=0,\; (V_2-12)d\xi_2=0,\\
&& (V_1-1)d\mu_1=0,\; (V_2-3)d\mu_2=0,\\
&& d\mu_1\geq 0,\; d\mu_2\geq 0,\; \mu_1(0)=0,\; \mu_2(0)=0,\\
&& s_1(p_1-p_2)=0,\\
&& s_2p_2=0,\\
&&p_1(0)=p_1(16),\quad p_2(0)=p_2(16),\\
&&\max_{w\in U} (-p_1w_1+p_2(-w_2+w_1)\\
&&+\lambda c (w_1(V_1+13-V_2)+w_2V_2)),\\
&& = -p_1u_1+p_2(-u_2+u_1)\\
&&+\lambda c (u_1(V_1+13-V_2)+u_2V_2),\\
&& \lambda+|p_1(16)|+|p_2(16)|+\mu_1(16)+\mu_2(16)=1.
\end{eqnarray*} 

Observe that $\dot{V}_1(t)\geq 0$. If $s_1(t)=0$, then we have $\dot{V}_1(t)>0$,  contradicting the periodicity of $V_1$. Hence, we have $s_1(t)>0$, $t\in [0,16]$. This implies that $V_1(t)\equiv 5$, $u_1(t)+s_1(t)\equiv 2$,  $p_1(t)\equiv p_2(t)$, and $d\mu_1=0$,  $t\in [0,16]$.

Suppose that $\lambda=0$. Then we have
\begin{eqnarray*}
&& dp_2=d\xi_2-d\mu_2,\\
&& s_2p_2=0,\\
&& p_2(0)=p_2(16),\\
&& \max_{w\in U} (-p_2w_2)=-p_2u_2,\\
&& |p_2(16)|+ \mu_2(16)>0.
\end{eqnarray*} 
 We claim that this is impossible. First, let us show that $V_2$ always reaches the  boundary of the interval $[3,12]$. Indeed, if $3<V_2(t)<12$, $t\in [0,16]$, then we have $d\xi_2=d\mu_2=0$ and $p_2(t)\equiv p_2(0)$. Thus, if $p_2(0)<0$, we have $\dot{V}_2(t)=-1$, which  implies that $V_2(\theta)=3$ for some $\theta< 16$. On the other hand, if $p_2(0)>0$, we have $\dot{V}_2(t)=3$, leading to the existence of some $\theta< 16$ such that $V_2(\theta)=12$. Hence, we have $p_2(16)=p_2(0)=0$,  contradicting the nontriviality condition.

Having proved that $V_2$ always reaches the boundary when $\lambda=0$, we next analyze what happens when such situations occur.  

Observe that the periodic solution $V_2\equiv 3$ is impossible. Indeed, in this case, we have $s_2\equiv 0$, $u_2\equiv 2$ and, as a consequence,  $p_2\equiv 0$. Hence, since $d\xi_2=0$, we obtain $d\mu_2=0$. This contradicts the nontriviality condition.

Now,  suppose that $V_2(\tau)>3$, for some $\tau\in [0, 16[$. Let us show that $V_2$ never reaches the value $3$. If there exists a first moment of time  $t=\theta\in ]\tau, 16]$ such that $V_2(\theta)=3$, then  $u_2(t)>2$, $t\in ]\theta-\delta,\theta [$, for some $\delta>0$.  This implies that  $p_2(t)\leq 0$,  $t\in ]\theta-\delta,\theta [$. Moreover, we have $d\mu_2(t)=0$, $t\in [\tau,\theta[$.  Since the case $V_2\equiv 3$ has already been excluded, we know that, once the lower bound of the state constraint is reached, $V_2$ must increase to preserve the periodicity of $V_2$. So,  there should exist a first instant $\theta_1$, such that  $\theta\leq \theta_1< 16$ and  $\dot V_2(t)>0$ for $t\in ]\theta_1,\theta_1+ \delta_1 [$. This means that $u_2(t)<2$, $t\in ]\theta_1,\theta_1+ \delta_1 [$ and,  it follows  from   $\max_{w\in U} (-p_2w_2)=-p_2u_2$, that  $p_2(t)\geq 0$,$t\in ]\theta_1,\theta_1+ \delta_1 [$. But this   is impossible  
since $d\mu_2\geq 0$ and $d\xi_2\cdot d\mu_2=0$. Thus, we must have $p_2(t)=0$ and $d\mu_2(t)=0$ for $t\in [\theta,16]$ and, consequently,   
  $|p_2(16)|+\mu_2(16)=0$, which is a contradiction. 

Let us show that the periodic solution $V_2\equiv 12$ is impossible. This implies that $p_2\geq 0$ to guarantee the condition $\dot{V}_2\geq 0$. The inequality $p_2>0$ is not possible, since the complementary condition $s_2p_2=0$ implies $s_2=0$. Hence $p_2\equiv 0$. This, together with the equality $\mu_2\equiv 0$, contradicts the nontriviality condition.

Finally, suppose that $V_2(\tau)<12$ for some $\tau\in [0,16[$. Let us show that $V_2$ never reaches the value $12$. We already know that $V_2>3$ and, therefore, $d\mu_2=0$. If there exists a first moment of time  $\theta\in]\tau, 16]$ such that $V_2(\theta)=12$, then  $u_2(t)\leq 0$, $t\in ]\theta-\delta,\theta [$, for some $\delta>0$.  This implies that  $p_2(t)\geq 0$,  $t\in ]\theta-\delta,\theta [$. As above, from the complementary condition $s_2p_2=0$, we see that $p_2(\tau_1)\leq 0$, for some $\tau_1\geq \theta$, to guarantee the periodicity of $V_2$. But, since $p_2$ is also periodic, we obtain $p_2\equiv 0$, contradicting the nontriviality condition. It follows from the above considerations that $\lambda>0$. 

Without loss of generality we set $\lambda=1$. In order to simplify the analysis of necessary conditions of optimality we introduce new adjoint variables 
\begin{eqnarray*}
&& Q_1=c(V_1+13-V_2),\\
&& Q_2=-p_2+cV_2.
\end{eqnarray*} 
In terms of these adjoint variables the necessary conditions take the form
\begin{eqnarray*}
&& dQ_1=c(u_2+s_2-2)dt +(V_1+13-V_2)dc,\\
&& dQ_2=-d\xi_2+c(u_2-u_1)dt+d\mu_2\\
&&+c(-u_2-s_2+2)dt+V_2dc,\\
&& \max_{w\in U} (Q_1w_1+Q_2w_2)=Q_1u_1+Q_2u_2.
\end{eqnarray*}
We highlight  the simplicity of the maximum principle condition in these new adjoint variables. 
Since $Q_1=c(V_1+13-V_2)>0$, we have $u_1\equiv 1$ and $s_1\equiv 1$. Thus, we have to analyze the following system:
\begin{eqnarray*}
&& \dot{V}_2=-u_2-s_2+2,\\
&& dQ_2=-d\xi_2+c dt+d\mu_2-cs_2dt+V_2dc,\\
&& \max_{w\in U} (Q_2w_2)= Q_2u_2.
\end{eqnarray*}
If $Q_2(\theta)>0$ and $V_2(\theta)<12$, then we have
$$
dQ_2=cdt+d\mu_2+V_2dc>0,\quad t\neq 16,
$$
 $u_2(t)=3$ and $\dot{V}_2(t)=-1-s_2=-1$, $t\in [\theta,16[$. Hence $Q_2(t)>0$, $t\in [\theta,16[$. 
 
Suppose that $V_2(t)=12$, $t\in [\theta_0,\theta_1]$. Then we have $u_2(t)+s_2(t)=2$,  $t\in [\theta_0,\theta_1]$. Hence $u_2(t)\leq 2$. The control $u_2$ appears in the functional with the positive coefficient $cV_2$. If $u_2<2$, then $\dot{V}_2>0$ and necessarily we must have $s_2>0$, which from the optimality of the process is senseless. We then conclude that  $u_2(t)=2$, $t\in [\theta_0,\theta_1]$. Hence, $s_2(t)=0$ and $Q_2(t)=0$, $t\in [\theta_0,\theta_1]$. Moreover,
$$
d\xi_2=cdt+V_2dt.
$$
  
The above analysis implies that, in general,  the variable $Q_2$ and the optimal control $u_2$ have the following forms:
$$
Q_2(t)=\left\{
\begin{array}{cl}
<0, & t\in [0,\tau_1[,\\
=0, & t\in [\tau_1,\tau_2],\\
>0, & t\in ]\tau_2,16[,\\
Q_2(0), & t=16.
\end{array}
\right.
$$
and
$$
u_2(t)=\left\{
\begin{array}{rl}
-1 & t\in [0,\tau_1[,\\
 2, & t\in [\tau_1,\tau_2],\\
3, & t\in ]\tau_2,16[.
\end{array}
\right.
$$
The jumps of $Q_2$ at $t=\tau_1$ and $t=\tau_2$, $dQ_2(\tau_1)$ and $dQ_2(\tau_2)$ are nonnegative.

Let $V_2(0)=V_2(16)=V^0<12$. We have
$$
12=V_2(\tau_1)=V^0+3\tau_1,
$$
$$
 V^0=V_2(16)=12-(16-\tau_2)=\tau_2-4.
$$
Hence $\tau_1=\frac{12-V^0}{3}\leq 3<6$ and $\tau_2=V^0+4\geq 7>6$.

Thus, we have to maximize the following quadratic function of $V^0$:
$$
\int_0^{16} c(u_1(V_1+13-V_2)+u_2V_2)dt
$$
$$
=\int_0^{\frac{12-V^0}{3}}3(5+13-2(V^0+3t))dt
+\int_{\frac{12-V^0}{3}}^63(6+2\cdot 12)dt
$$
$$
+\int_6^{V^0+4}11(6+2\cdot 12)dt
$$
$$
+\int_{V^0+4}^{16} 11((5+13-(12-(t-(V^0+4)))
$$
$$
+3(12-(t-(V^0+4))))dt.
$$
It attains its maximum at $V^0=\frac{36}{5}$. The corresponding $Q_2$ has the form
$$
Q_2(t)=\left\{
\begin{array}{cl}
-(12-V^0)+3t,& t\in [0,\tau_1],\\
0, & t\in [\tau_1,\tau_2],\\
11(t-(V^0+4)),  & t\in [\tau_2,16[,\\
-(12-V^0), & t=16,
\end{array}
\right.
$$
where $\tau_1=\frac{12-V^0}{3}$ and $\tau_2=V^0+4$. The spillway $s_2(t)$ is always zero.

\section{Conclusions}

 We have derived necessary conditions of optimality for the optimal control problem arising in optimization of 
  the   power production profit for cascades of hydro-electric power plants with reversible turbines and uncontrolled spillways. 
    We have illustrated the effectiveness of our necessary conditions solving analytically an example involving two hydro-electric power plants. 

A remarkable feature of the proof of our main result is that the solution of the original problems is  approximated by the solutions of a family of standard optimal control problems. For these latter problems,  there is a substantial literature on numerical methods, and, therefore, our approximation techniques may be used to construct numerical methods to solve real-world problems for cascades with several hydro-electric power plants and uncontrolled spillways.

\section*{Appendix}

{\em Proof of Lemma \ref{APP}}. We prove the lemma by induction on the number of equations. 
The first equation in \eqref{A1} has the form
\begin{equation}
\label{Y1}
\dot{V}_1^{\gamma}=A_1-u_1^{\gamma}-s^{\gamma}_1.
\end{equation}
If system \eqref{A1} were reduced to the above equation, the result would  follow from a more general result \cite[Theorem 4]{S}. Since this is not the case, and to make the presentation self-contained, we  give a detailed simple proof for equation \eqref{Y1}.
 
Multiplying  \eqref{Y1} by $\dot{V}_1^{\gamma}$ we get
$$
|\dot{V}_1^{\gamma}|^2=(A_1-u_1^{\gamma})\dot{V}_1^{\gamma}-\frac{d}{dt} e^{\gamma (V_1^{\gamma}-V_1^M)}.
$$
Integrating this equality and using the inequality $ab\leq\frac{a^2+b^2}{2}$, we obtain
$$
\int_0^T |\dot{V}_1^{\gamma}|^2 dt\leq \frac{T(A_1+\max\{|u_1^m|,|u_1^M|\} )^2}{2}
$$
$$
+\frac{1}{2}\int_0^T |\dot{V}_1^{\gamma}|^2 dt-
e^{\gamma (V_1^{\gamma}(T)-V_1^M)}+e^{\gamma (V_1^{\gamma}(0)-V_1^M)}.
$$
Hence, we have
\begin{equation}
\label{Y2}
\int_0^T |\dot{V}_1^{\gamma}|^2 dt\leq T(A_1+\max\{|u_1^m|,|u_1^M|\} )^2 +2.
\end{equation}
Integrating \eqref{Y1} we obtain
$$
|V_1^{\gamma}(t_2)-V_1^{\gamma}(t_1)|\leq \int_{t_1}^{t_2}|\dot{V}_1^{\gamma}|dt
$$
$$
\leq \sqrt{t_2-t_1}\left( \int_{t_1}^{t_2}|\dot{V}_1^{\gamma}|^2 dt\right)^{\frac{1}{2}}.
$$
This, together with \eqref{Y2}, allows us to apply the Arzela-Ascoli theorem. Hence, without loss of generality, the sequence  $V_1^{\gamma}$ uniformly tends to $\tilde{V}_1\in C([0,T],R)$.

It is an easy matter to see that  \eqref{Y1} and \eqref{Y2} imply that the sequence of functions $s_1^{\gamma}$ is bounded in $L_2([0,T],R)$.
Therefore, without loss of generality, the sequence $s_1^{\gamma}$ converge in weak topology of $L_2([0,T],R)$ to some function $\tilde{s}_1$.

Integrating \eqref{Y1} again, we also get
$$
V_1^{\gamma}(t)=V_1^{\gamma}(0)+\int_0^t (A_1-u_1^{\gamma}-s_1^{\gamma})dt.
$$
Passing to the limit in this equality we obtain
$$
\tilde{V}_1(t)=V_1^0+\int_0^t (A_1-\tilde{u}_1-\tilde{s}_1)dt.
$$
 Thus, we have
 \begin{equation}
 \label{Y3}
 \dot{\tilde{V}}_1=A_1-\tilde{u}_1-\tilde{s}_1.
 \end{equation}
 
 Let us show that $\title{V}_1(t)\leq V_1^M$. Indeed, if $V_1^{\gamma}(t_*)=V_1^M$ and the derivative $\dot{V}_1^{\gamma}(t_*)$ exists, then, since $s_1^{\gamma}(t_*)=\gamma$, we have
 $$
 \dot{V}_1^{\gamma}(t_*)=A_1-u_1^{\gamma}-s_1^{\gamma}=A_1-u_1^{\gamma}-\gamma<0
 $$
 for $\gamma >  A_1+\max\{|u_1^m|,|u_1^M|\}$. Hence, $V_1^{\gamma}(t)\leq V_M$ for such $\gamma$. This implies that $\tilde{V}_1(t)\leq V_1^M$. 
 
 Now, let us show that $\tilde{s}_1(\cdot)\in L_{\infty}([0,T],R)$. If $\tilde{V}_1(t)<V_1^M$, then we obviously have $\tilde{s}_1(t)=0$. Assume that $\tilde{V}_1(t_*)=V_1^M$ and the derivative $\dot{\tilde{V}}_1(t_*)$ exists. Then, since $\tilde{V}_1(t)\leq V_1^M$, for all $t\in [0,T]$, we have
 $$
 0=\dot{\tilde{V}}_1(t_*)=A_1-\tilde{u}_1(t_*)-\tilde{s}_1(t_*),
 $$
 proving our claim since $|\tilde{s}_1(t_*)|\leq  A_1+\max\{|u_1^m|,|u_1^M|\}$.

To prove the uniqueness, observe that  if there exist two admissible solutions $V_1$, $s_1$, $V'_1$, $s'_1$ to equation \eqref{Y2}, then we have
 $$
\dot{V}_1-\dot{V}'_1=-s_1+s_1'.
$$
Multiplying this by $V_1-V_1'$, we get
$$
\frac{1}{2}\frac{d}{dt}(V_1-V_1')^2=(-s_1+s_1')(V_1-V_1')
$$
$$
=(-s_1+s_1')(V_1-V_1^M-V_1'+V_1^M)
$$
$$
=s_1(V_1'-V_1^M)+s_1'(V_1-V_1^M)\leq 0.
$$
Therefore, $V_1(t)$ and $s_1(t)$ are uniquely defined.

\vspace{5mm}

Since the inclusion $j\in J(i)$ implies that $j<i$, arguing as above, we prove, by induction,  that,  
without loss of generality, when $\gamma\to\infty$,  the functions $s^{\gamma}_i$ converge in the weak topology of $L_2([0,T],R)$ to some essentially bounded measurable functions $\tilde{s}_i$ and the functions $V_i^{\gamma}$ uniformly tend to $\tilde{V}_i\leq V_i^M$ satisfying
$$
\dot{\tilde{V}}_i=A_i-\tilde{u}_i-\tilde{s}_i+\sum_{j\in J(i)}\left(\tilde{u}_j+\tilde{s}_j\right),
$$
where $\tilde{s}_i\geq 0$ and $\tilde{s}_i(\tilde{V}_i-V_i^M)=0$.  In the proof of the inequality $\tilde{V}_i\leq V_i^M$,  it is important to notice that, when $V_i^{\gamma}(t)=V_i^M$, then, for sufficiently large $\gamma$, the inequality
$$
\dot{V}_i^{\gamma}=A_i-u_i^{\gamma}-s^{\gamma}_i +\sum_{j\in J(i)}\left( u_j^{\gamma}+s^{\gamma}_j\right)
$$
$$
\leq A_i+\max\{|u_i^m|,|u_i^M|\}-\gamma^i 
$$
$$
+ \sum_{j\in J(i)}\left(\max\{|u_j^m|,|u_j^M|\}+\gamma^j\right)<0
$$
holds.

\vspace{5mm}

Now, let us show that the Cauchy problem for system  \eqref{s1} -- \eqref{s2} has a unique solution. Indeed, let $u$ be an admissible control and let $V$, $s$, $V'$, $s'$ be two admissible solutions to problem \eqref{s1} and \eqref{s2} corresponding to the same $u$ and with the same initial condition. We have
\begin{eqnarray*}
&&  \dot{V}=A-u-s+{\cal M}(u+s),\\
&& V\leq V^M,\;\; s\geq 0,\;\; {\rm diag}(s)(V-V^M)=0,\\
&& V(0)=V^0,
\end{eqnarray*}
and
\begin{eqnarray*}
&&  \dot{V}'=A-u-s'+{\cal M}(u+s'),\\
&& V'\leq V^M,\;\; s'\geq 0,\;\; {\rm diag}(s')(V'-V^M)=0,\\
&& V'(0)=V^0.
\end{eqnarray*}
Suppose that we have proved the uniqueness of the first $i$ components of the vectors $V$ and $s$.  Then (see the definition of $J(i)$) we have
$$
\dot{V}_{i+1}-\dot{V}'_{i+1}=-s_{i+1}+s_{i+1}'.
$$
Arguing exactly as in the case $i=1$, we see that $V_{i+1}(t)$ and $s_{i+1}(t)$ are uniquely defined. Therefore, the solution to problem \eqref{s1} and \eqref{s2} is unique.
$\Box$


\bigskip

{\em Proof of Lemma \ref{lem1}}.
We fix $\alpha$ and $\epsilon$. 
Denote by $\hat{V}^{\gamma}$ the solution to \eqref{ap2} corresponding to the control $\hat{u}$ with the initial condition $\hat{V}(0)$. By Lemma \ref{APP}, the sequence $\hat{V}^{\gamma}$ converges to $\hat{V}$ as $\gamma\to\infty$. Hence, $\hat{V}^{\gamma}$ is an admissible solution to problem \eqref{ap1} -- \eqref{ap4} for large $\gamma$.  Therefore, applying  classical existence theorems (see, e.g., Chapter III, Theorem 2.1 and Corollary 4.1 in \cite{FR}) we see that there exists a global optimal control process  $(\tilde{u}^{\gamma,\epsilon,\alpha},\tilde{V}^{\gamma,\epsilon,\alpha},\tilde{s}^{\gamma,\epsilon,\alpha})$.  Without loss of generality the sequence $\tilde{u}^{\gamma,\epsilon,\alpha}$ converges in the weak-* topology of $L_{\infty}([0,T],R^{\cal I})$ to some control $\tilde{u}^{\epsilon,\alpha}\in  [u^m,u^M]$,  the sequence $\tilde{V}^{\gamma,\epsilon,\alpha}$ uniformly converges to $\tilde{V}^{\epsilon,\alpha}$ as $\gamma\to\infty$ and the sequence $\tilde{s}^{\gamma,\epsilon,\alpha}$ converges in the weak topology of $L_2([0,T],R^{\cal I})$ to some function $\tilde{s}^{\epsilon,\alpha}\in L_{\infty}([0,T],R^{\cal I})$.
 
We have
$$
{\cal J}(\tilde{u}^{\gamma,\epsilon,\alpha},\tilde{V}^{\gamma,\epsilon,\alpha})\leq {\cal J}(\hat{u},\hat{V}^{\gamma})
$$
$$
=
 - \int_0^T\sum_{j=1}^{{\cal I}} c \hat{u}_j\left(\frac{\hat{V}_j^{\gamma}}{S_j}+h_j-\frac{\hat{V}_{J^{-1}(j)}^{\gamma}}{S_{J^{-1}(j)}}-h_{J^{-1}(j)}\right) dt
$$
$$
+ \int_0^T\sum_{j=1}^{{\cal I}}\frac{e^{-\gamma (\hat{V}_j^{\gamma}-V_j^m+\epsilon)}-1}{\sqrt{\gamma}} dt
$$
From this and \eqref{Vbound}, for large $\gamma$, we get
$$
 \int_0^T\sum_{j=1}^{{\cal I}}\frac{e^{-\gamma (\tilde{V}_j^{\gamma,\epsilon,\alpha}-V_j^m+\epsilon)}-1}{\sqrt{\gamma}}dt
 $$
 $$
\leq
 - \int_0^T\sum_{j=1}^{{\cal I}} c \hat{u}_j\left(\frac{\hat{V}_j^{\gamma}}{S_j}+h_j-\frac{\hat{V}_{J^{-1}(j)}^{\gamma}}{S_{J^{-1}(j)}}-h_{J^{-1}(j)}\right) dt
$$
$$
+\int_0^T\sum_{j=1}^{{\cal I}} c \tilde{u}_j^{\gamma,\epsilon,\alpha}\left(\frac{\tilde{V}_j^{\gamma,\epsilon,\alpha}}{S_j}+h_j-\frac{\tilde{V}_{J^{-1}(j)}^{\gamma,\epsilon,\alpha}}{S_{J^{-1}(j)}}-h_{J^{-1}(j)}\right) dt
$$
$$
\leq M_7.
$$
Let us show that $\tilde{V}_i^{\epsilon,\alpha}(t)\geq V_i^m-\epsilon$, $i=1,\ldots,{\cal I}$. Set  $ \Delta=\{ t\in [0,T] :~ \tilde{V}_i^{\epsilon,\alpha}(t)<V_i^m-\epsilon\}$.  Passing to the limit as $\gamma\to\infty$ in the inequality
$$
\int_{\Delta} (V_i^m-\epsilon-\tilde{V}_i^{\gamma,\epsilon,\alpha})dt-\frac{T({\cal I}-1)}{\gamma} 
$$
$$
\leq \int_0^T \sum_{j\neq i}\frac{e^{-\gamma (\tilde{V}_j^{\gamma,\epsilon,\alpha}-V_j^m+\epsilon)}-1}{\gamma} dt \leq \frac{M_8}{\sqrt{\gamma}},
$$
we get
$$
0<\int_{\Delta} (V_i^m-\epsilon-\tilde{V}_i^{\epsilon,\alpha})dt\leq 0,
$$
a contradiction. Hence, $\tilde{V}_i^{\epsilon,\alpha}\geq V_i^m-\epsilon$, $i=1,\ldots,{\cal I}$. Thus, $(\tilde{u}^{\epsilon,\alpha},\tilde{V}^{\epsilon,\alpha},\tilde{s}^{\epsilon,\alpha})$ is an admissible solution to problem $({\cal P}_{\epsilon,\alpha})$. $\Box$

\bigskip

{\em Proof of Lemma \ref{lem2}}. 
Observe that if $\tilde{V}_1^{\epsilon,\alpha}(t)=V_1^M$ and the derivative $\dot{\tilde{V}}_1^{\epsilon,\alpha}(t)$ exists, then $\dot{\tilde{V}}_1^{\epsilon,\alpha}(t)=0$. Hence, we have $\tilde{s}_1^{\epsilon,\alpha}=A_1-\tilde{u}_1^{\epsilon,\alpha}$. By induction, we get  
$$
\tilde{s}_i^{\epsilon,\alpha}=A_i-\tilde{u}_i^{\epsilon,\alpha}+\sum_{j\in J(i)}( \tilde{u}_j^{\epsilon,\alpha}+\tilde{s}_j^{\epsilon,\alpha}),
$$
whenever $\tilde{V}_i^{\epsilon,\alpha}(t)=V_i^M$ and the derivative $\dot{\tilde{V}}_i^{\epsilon,\alpha}(t)$ exists. Thus, we see that the sequences $\tilde{s}_i^{\epsilon,\alpha}\in L_{\infty}$ are bounded. 

Since
$$
\tilde{V}_i^{\epsilon,\alpha}(t)=\tilde{V}_i^{\epsilon,\alpha}(0)+\int_0^t\left( A_i-\tilde{u}_i^{\epsilon,\alpha}-\tilde{s}_i^{\epsilon,\alpha}
+\sum_{j\in J(i)}( \tilde{u}_j^{\epsilon,\alpha}+\tilde{s}_j^{\epsilon,\alpha})\right) dt,
$$
it is a simple matter to see that, without loss of generality, it converges to an admissible solution  $(\tilde{u}^{\alpha},\tilde{V}^{\alpha},\tilde{s}^{\alpha})$ of  problem \eqref{P1} -- \eqref{P5}, as $\epsilon\to 0$.

Let us go back to the sequence $(\tilde{u}^{\gamma,\epsilon,\alpha},\tilde{V}^{\gamma,\epsilon,\alpha})$ of the global optimal solutions to problem \eqref{ap1} -- \eqref{ap4}.  From the inequality
$$
{\cal J}(\tilde{u}^{\gamma,\epsilon,\alpha},\tilde{V}^{\gamma,\epsilon,\alpha})\leq {\cal J}(\hat{u},\hat{V}^{\gamma})
$$
we have
$$
\frac{1}{2}|\tilde{V}^{\gamma,\epsilon,\alpha}(0)-\hat{V}(0)|^2 
$$
$$
- \int_0^T\sum_{j=1}^{{\cal I}} c \tilde{u}_j^{\gamma,\epsilon,\alpha}\left(\frac{\tilde{V}_j^{\gamma,\epsilon,\alpha}}{S_j}+h_j-\frac{\tilde{V}_{J^{-1}(j)}^{\gamma,\epsilon,\alpha}}{S_{J^{-1}(j)}}-h_{J^{-1}(j)}\right) dt
$$
$$
 - \frac{T{\cal I}}{\sqrt{\gamma}}+\int_0^T\alpha |\tilde{u}^{\gamma,\epsilon,\alpha}-\hat{u}|^2 dt
 $$
 $$
\leq
 - \int_0^T\sum_{j=1}^{{\cal I}} c \hat{u}_j\left(\frac{\hat{V}_j^{\gamma}}{S_j}+h_j-\frac{\hat{V}_{J^{-1}(j)}^{\gamma}}{S_{J^{-1}(j)}}-h_{J^{-1}(j)}\right) dt,
$$
whenever $\gamma$ is sufficiently large. Passing to the limit when $\gamma$ tends to infinity, we get
$$
\frac{1}{2}|\tilde{V}^{\epsilon,\alpha}(0)-\hat{V}(0)|^2 
$$
$$
- \int_0^T\sum_{j=1}^{{\cal I}} c \tilde{u}_j^{\epsilon,\alpha}\left(\frac{\tilde{V}_j^{\epsilon,\alpha}}{S_j}+h_j-\frac{\tilde{V}_{J^{-1}(j)}^{\epsilon,\alpha}}{S_{J^{-1}(j)}}-h_{J^{-1}(j)}\right) dt
$$
$$
+\int_0^T\alpha |\tilde{u}^{\epsilon,\alpha}-\hat{u}|^2 dt
=\frac{1}{2}|\tilde{V}^{\epsilon,\alpha}(0)-\hat{V}(0)|^2 
$$
$$
- \int_0^T\sum_{j=1}^{{\cal I}} c \tilde{u}_j^{\epsilon,\alpha}\left(\frac{(\tilde{V}_j^{\epsilon,\alpha}-\tilde{V}_j^{\alpha})}{S_j}+h_j
 -\frac{(\tilde{V}_{J^{-1}(j)}^{\epsilon,\alpha}-\tilde{V}_{J^{-1}(j)}^{\alpha})}{S_{J^{-1}(j)}}-h_{J^{-1}(j)}\right) dt
$$
$$
- \int_0^T\sum_{j=1}^{{\cal I}} c \tilde{u}_j^{\epsilon,\alpha}\left(\frac{\tilde{V}_j^{\alpha}}{S_j}+h_j-\frac{\tilde{V}_{J^{-1}(j)}^{\alpha}}{S_{J^{-1}(j)}}-h_{J^{-1}(j)}\right) dt
$$
$$
+\int_0^T\alpha |\tilde{u}^{\epsilon,\alpha}-\hat{u}|^2 dt
$$
$$
\leq 
 - \int_0^T\sum_{j=1}^{{\cal I}} c \hat{u}_j\left(\frac{\hat{V}_j}{S_j}+h_j-\frac{\hat{V}_{J^{-1}(j)}}{S_{J^{-1}(j)}}-h_{J^{-1}(j)}\right) dt,
$$
Since the convergence of $\tilde{u}^{\epsilon,\alpha}$ in the weak-* topology of $L_{\infty}$ implies its convergence in the weak topology of $L_2$ and the $L_2$-norm is lower semi-continuous with respect to weak convergence, passing to the limit when $\epsilon$ tends to zero, we get
$$
\frac{1}{2}|\tilde{V}(0)-\hat{V}(0)|^2
$$
$$
 - \int_0^T\sum_{j=1}^{{\cal I}} c \tilde{u}_j\left(\frac{\tilde{V}_j^{\alpha}}{S_j}+h_j-\frac{\tilde{V}_{J^{-1}(j)}^{\alpha}}{S_{J^{-1}(j)}}-h_{J^{-1}(j)}\right) dt
$$
$$
+\int_0^T\alpha |\tilde{u}-\hat{u}|^2 dt
$$
$$
\leq 
 - \int_0^T\sum_{j=1}^{{\cal I}} c \hat{u}_j\left(\frac{\hat{V}_j}{S_j}+h_j-\frac{\hat{V}_{J^{-1}(j)}}{S_{J^{-1}(j)}}-h_{J^{-1}(j)}\right) dt
 $$
 $$
 \leq 
 - \int_0^T\sum_{j=1}^{{\cal I}} c \tilde{u}_j\left(\frac{\tilde{V}_j^{\alpha}}{S_j}+h_j-\frac{\tilde{V}_{J^{-1}(j)}^{\alpha}}{S_{J^{-1}(j)}}-h_{J^{-1}(j)}\right) dt.
$$
Hence $\tilde{u}=\hat{u}$ and $\tilde{V}=\hat{V}$. $\Box$

\bigskip

{\em Proof of Theorem \ref{th_main}}. 
Now, we are in a position to prove our main result. 
 Let $(u,V,s)$ be a global optimal solution to problem  \eqref{ap1} -- \eqref{ap4}. We apply the classical Pontryagin maximum principle to problem \eqref{ap1} -- \eqref{ap4}. To avoid complicated notations, we omit the indexes $\gamma$, $\epsilon$ and $\alpha$ in the optimal control process, the adjoint variables and multipliers. 
 
A word of caution while reading the proof of this theorem. Omitting the indexes in a sequence $F^{(\gamma,\epsilon,\alpha)}$ comes with a price for, when taking limits, we use the same letter $F$ to denote the limit. For example, instead of writing $\lim_{\gamma\to\infty} F^{(\gamma,\epsilon,\alpha)}=F^{(\epsilon,\alpha)}$, we write
 $\lim_{\gamma\to\infty} F=F$, etc..
 
 Thus, the Hamiltonian is 
$$
H=\sum_{i=1}^{\cal I} q_i\left( A_i-w_i-\gamma^i e^{\gamma^i(V_i-V_i^M)}
+\sum_{j\in J(i)}\left( w_j+\gamma^j e^{\gamma^j(V_j-V_j^M)}\right) \right)
$$
$$
+\lambda^0 \left( \sum_{j=1}^{\cal I} cw_j\left( \frac{V_j}{S_j}+h_j - \frac{V_{J^{-1}(j)}}{S_{J^{-1}(j)}}-h_{J^{-1}(j)}\right)\right.
$$
$$
\left. -
\sum_{j=1}^{\cal I} \frac{e^{-\gamma (V_j-V_j^m+\epsilon)}-1}{\sqrt{\gamma}}-\alpha|w-\hat{u}|^2\right).
$$ 
Here, we denote a control by $w$ to distinguish it from the optimal control $u$. 
The corresponding adjoint system reads
$$
\dot{q}_i=\gamma^{2i} e^{\gamma^i(V_i-V_i^M)}(q_i-q_{J^{-1}(i)})
$$
$$
+\lambda^0\left( -c\;\frac{u_i-\sum_{j\in J(i)}u_j}{S_i} 
-\sqrt{\gamma}e^{-\gamma (V_i-V_i^m+\epsilon)}\right).
$$
Moreover, we have the transversality conditions
$$
q_i(0)=\lambda^0(V_i(0)-\hat{V}_i(0))+\nu^0(V_i(0)-V_i(T)),
$$
$$
 q_i(T)=\nu^0(V_i(0)-V_i(T)),
$$
 the complementary condition
$$
\nu^0(|V(0)-V(T)|^2-2\epsilon)=0,\quad \nu^0\geq 0,
$$
and the nontriviality condition
$$
\lambda^0+\nu^0>0.
$$
The maximum principle reads: almost everywhere the optimal control $u(t)$ solves the maximization problem
$$
\sum_{i=1}^{\cal I} q_i\left( -w_i+\sum_{j\in J(i)}w_j \right)
$$
$$
+\lambda^0 \left( \sum_{j=1}^{\cal I} cw_j\left( \frac{V_j}{S_j}+h_j
 - \frac{V_{J^{-1}(j)}}{S_{J^{-1}(j)}}-h_{J^{-1}(j)}\right) -\alpha|w-\hat{u}|^2\right)\to\max_{w\in U}.
$$ 

\vspace{5mm}

Before passing to the limit we normalize the adjoint variables and the multipliers $\lambda^0$. 
Observe that 
$$
\lambda^0+\nu^0|V(0)-V(T)| +\lambda^0\sum_{i=1}^{\cal I} \int_0^T\sqrt{\gamma}e^{-\gamma (V_i-V_i^m+\epsilon)}dt>0.
$$
Indeed, if this is not the case, then we have $\lambda^0=0$ and $|V(0)-V(T)|=0$. Then from the complementary condition we obtain $\nu^0=0$, a contradiction.

Denote  the following: 
$$
\frac{1}{ \lambda^0+\nu^0|V(0)-V(T)| +\lambda^0\sum_{i=1}^{\cal I} \int_0^T\sqrt{\gamma}e^{-\gamma (V_i-V_i^m+\epsilon)}dt}
$$
by $\sigma$. Set
$$
\lambda=\sigma\lambda^0,\quad p=\sigma q,\quad \mu_i=\sigma\lambda^0 \int_0^t\sqrt{\gamma}e^{-\gamma (V_i-V_i^m+\epsilon)}d\tau.
$$
Hence we get
$$
\dot{p}_i=\gamma^{2i} e^{\gamma^i(V_i-V_i^M)}(p_i-p_{J^{-1}(i)})-\lambda c\frac{u_i-\sum_{j\in J(i)}u_j}{S_i} -\dot{\mu}_i,
$$
$$
p(0)=\lambda (V(0)-\hat{V}(0))+p(T).
$$
$$
\lambda+|p(T)|+\sum_{i=1}^{\cal I} \mu_i(T)=1.
$$
Moreover, almost everywhere the optimal control $u(t)$ solves the maximization problem
$$
\sum_{i=1}^{\cal I} p_i\left( -w_i+\sum_{j\in J(i)}w_j \right)
$$
$$
+\lambda  \sum_{j=1}^{\cal I} cw_j\left( \frac{V_j}{S_j}+h_j - \frac{V_{J^{-1}(j)}}{S_{J^{-1}(j)}}-h_{J^{-1}(j)}\right)
$$
$$
-\lambda\alpha|w-\hat{u}|^2\to\max_{w\in U}.
$$ 

\vspace{5mm}

We fix $\epsilon$ and $\alpha$. Recall that $V$, $u$, $p$ and $\mu$ depend on $\gamma$ but we omit this index.  Now, we are in a position to pass to the limit in these normalized necessary conditions, as $\gamma\to\infty$.  
This is done by induction on the number of equations, as in the proof of Lemma \ref{APP}, but, this time, backwards, i.e., from the last equation to the first one.  

For the adjoint variable corresponding to $i={\cal I}$, we have (a case similar to \cite[Theorem 5]{S}):
\begin{equation}
\label{ppI}
\dot{p}_{\cal I}=\gamma^{2{\cal I}} e^{\gamma^{\cal I}(V_{\cal I}-V_{\cal I}^M)}p_{\cal I} -\lambda c\frac{u_{\cal I}-\sum_{j\in J({\cal I})}u_j}{S_{\cal I}}-\dot{\mu}_{\cal I}.
\end{equation}
Integrating \eqref{ppI}  we get
$$
\frac{d}{dt}|p_{\cal I}| =\gamma^{2{\cal I}} e^{\gamma^{\cal I}(V_{\cal I}-V_{\cal I}^M)}|p_{\cal I}|
$$
$$
 -\left(\lambda c\frac{u_{\cal I}-\sum_{j\in J({\cal I})}u_j}{S_{\cal I}}+\dot{\mu}_{\cal I}\right) {\rm sign}(p_{\cal I})
$$
$$
\geq -\left(\lambda c\frac{u_{\cal I}-\sum_{j\in J({\cal I})}u_j}{S_{\cal I}}+\dot{\mu}_{\cal I}\right) {\rm sign}(p_{\cal I}).
$$
From the Gronwall inequality we obtain
$$
|p_{\cal I}(t)|\leq |p_{\cal I}(T)|+M_9\leq M_9+1.
$$
Hence
$$
\int_0^T \gamma^{2{\cal I}} e^{\gamma {\cal I}(V_{\cal I}-V_{\cal I}^M)}|p_{\cal I}| dt=|p_{\cal I}(T)|-|p_{\cal I}(0)|
$$
$$
+\int_0^T\left( \lambda c\frac{u_{\cal I}-\sum_{j\in J({\cal I})}u_j}{S_{\cal I}}+\dot{\mu}_{\cal I} \right) {\rm sign}(p_{\cal I}) dt.
$$
Therefore
$$
\int_0^T \gamma^{2{\cal I}} e^{\gamma^{\cal I}(V_{\cal I}-V_{\cal I}^M)}|p_{\cal I}| dt\leq M_{10}
$$
and, hence,
$$
\int_0^T|\dot{p}_{\cal I}|dt\leq M_{11}.
$$
Thus, the sequences  $\mu_{\cal I}$, $\int_0^t \gamma^{2{\cal I}} e^{\gamma {\cal I}(V_{\cal I}-V_{\cal I}^M)} p_{\cal I} d\tau$ and $ \int_0^t \dot{p}_{\cal I}^{\gamma}d\tau$ are bounded in  $BV([0,T],R)$. Without loss of generality, they converge in the  weak-* topology to some functions $\mu_{\cal I}$, $\xi_{\cal I}$ and $p_{\cal I}\in BV([0,T],R)$,  the controls $u$ converge in the weak-* topology of $L_{\infty}([0,T],R^{\cal I})$ and the solutions $V$ converge in the norm of $C([0,T],R^{\cal I})$. Passing to the limit in \eqref{ppI} we obtain
$$
d{p}_{\cal I}=d\xi_{\cal I}-\lambda c\frac{u_{\cal I}-\sum_{j\in J({\cal I})}u_j}{S_{\cal I}}dt- d\mu_{\cal I},
$$
$$
 (V_{\cal I}-V_{\cal I}^M)d\xi_{\cal I}=0,\quad (V_{\cal I}-V_{\cal I}^m+\epsilon)d\mu_{\cal I}=0,\quad d\mu_{\cal I}\geq 0.
$$

Now, we continue the proof by induction over $i<{\cal I}$. Suppose that 
$$
\int_0^T|\dot{p}_{J^{-1}(i)}|dt\leq M_{12}.
$$
Recall that $J^{-1}(i)>i$. We have
$$
\frac{d}{dt}|p_i-p_{J^{-1}(i)}|=\gamma^{2i} e^{\gamma^i(V_i-V_i^M)}|p_i- p_{J^{-1}(i)}|
 $$
 $$
 -\left(\dot{p}_{J^{-1}(i)}
+\lambda c\frac{u_i-\sum_{j\in J(i)}u_j}{S_i}  +\dot{\mu}_i \right) {\rm sign}(p_i-p_{J^{-1}(i)})
$$
$$
\geq 
-\left( \dot{p}_{J^{-1}(i)}
+\lambda c\frac{u_i-\sum_{j\in J(i)}u_j}{S_i} +\dot{\mu}_i\right) {\rm sign}(p_i-p_{J^{-1}(i)}).
$$
Integrating this inequality  and using the Gronwall inequality we get
$$
|p_i(t)-p_{J^{-1}(i)}(t)|\leq |p_i(T)-p_{J^{-1}(i)}(T)|+M_{13}.
$$
Hence
\begin{equation}
\label{sp}
\int_0^T \gamma^{2i} e^{\gamma^i(V_i-V_i^M)}|p_i-p_{J^{-1}(i)}|dt\leq M_{14}
\end{equation}
and
$$
\int_0^T|\dot{p}_i|dt\leq M_{15}.
$$

Arguing as above and passing to the limit we see that there exist $u\in L_{\infty}([0,T],R^{\cal I})$, $V\in C([0,T],R^{\cal I})$, $\mu_i$, $\xi_i$ and $p_i\in BV([0,T],R)$ satisfying
$$
d{p}_i=d\xi_i-\lambda c \frac{u_i-\sum_{j\in J(i)}u_j}{S_i} dt- d\mu_i,
$$
$$
 (V_i-V_i^M)d\xi_i=0,\quad (V_i-V_i^m+\epsilon)d\mu_i=0,\quad d\mu_i\geq 0.
$$
Dividing \eqref{sp} by $\gamma$ and passing to the limit, we get
$$
s_i(p_i-p_{J^{-1}(i)})=0.
$$

All the constructed functions $\mu_i$,  $\xi_i$ and $p_i$,  $i=1,\ldots,{\cal I}$, are bounded in the space $BV([0,T],R)$. They  depend on $\epsilon$ and $\alpha$.

\vspace{5mm}

Now, we consider $\epsilon\downarrow 0$. Since the parameter $\epsilon$ appears only in the equalities  $(V_i-V_i^m+\epsilon)d\mu_i=0$, and the functions $\mu_i$,  $\xi_i$ and $p_i$,  $i=1,\ldots,{\cal I}$, constructed above, are bounded in the space $BV([0,T],R)$, without loss of generality they converge in the weak-* topology of the space $BV([0,T],R)$. Passing to the limit, we see that there exist functions $u\in L_{\infty}([0,T],R^{\cal I})$, $V\in C([0,T],R^{\cal I})$, $\mu_i$,  $\xi_i$ and $p_i$,  $i=1,\ldots,{\cal I}$, belonging to the space $BV([0,T],R)$ and satisfying
$$
d{p}_i=d\xi_i-\lambda c \frac{u_i-\sum_{j\in J(i)}u_j}{S_i} dt- d\mu_i,
$$
$$
 (V_i-V_i^M)d\xi_i=0,\quad (V_i-V_i^m)d\mu_i=0,\quad d\mu_i\geq 0.
$$
These functions depend on the parameter $\alpha$.


Next, using the boundness of the constructed functions $\mu_i$,  $\xi_i$ and $p_i$,  $i=1,\ldots,{\cal I}$,   we pass to the limit as $\alpha\downarrow 0$ as above. This ends the proof. $\Box$

\bibliographystyle{plain}
\bibliography{ref_Barragem}

\end{document}